\documentclass{elsart}

\usepackage{amsmath}
\usepackage{amssymb}
\journal{J. Combin. Theory Ser. A} 

\newcommand{\xv}{\mathbf{x}}

\newcommand{\av}{\mathbf{a}}

\newcommand{\bv}{\mathbf{b}}

\newcommand{\ev}{\mathbf{e}}
\newcommand{\zerov}{\mathbf{0}}
\newcommand{\s}{\sigma_n(\av)}
\newcommand{\qbin}[2]{\genfrac{[}{]}{0pt}{}{#1}{#2}}
\newcommand{\gp}[3]{\qbin{#1}{#2}_{#3}}

\newenvironment{DC}{\textbf{Dyson's conjecture} \em  }{}
\newenvironment{qDC}{\textbf{Andrews' $q$-Dyson conjecture} \em  }{}

\numberwithin{equation}{section}
\numberwithin{thm}{section}

\begin{document}
\begin{frontmatter}
\title{Disturbing the Dyson Conjecture,\\ in a \emph{Generally} GOOD Way}
\author{Andrew V. Sills}
\address{Department of Mathematics, Rutgers University,
         Hill Center, Busch Campus, Piscataway, NJ 08854-8019 USA}
\ead{asills@math.rutgers.edu}
\ead[url]{www.math.rutgers.edu/\~{}asills}
\date{4 November 2005}

\begin{abstract}
  Dyson's celebrated constant term conjecture 
({\em J. Math. Phys.}, 3 (1962): 140--156) states that the constant term
in the expansion of $\prod_{1\leqq i\neq j\leqq n} (1-x_i/x_j)^{a_j}$
is the multinomial coefficient $(a_1 + a_2 + \cdots + a_n)!/
(a_1! a_2! \cdots a_n!)$.
The definitive proof was given by I. J. Good
({\em J. Math. Phys.}, 11 (1970) 1884).
Later, Andrews extended Dyson's conjecture to a $q$-analog
({\em The Theory and Application of Special Functions},
(R. Askey, ed.), New York: Academic Press, 191--224,  1975.)
In this paper, closed form expressions are given for the coefficients of
several other terms in the Dyson product, and are proved using an
extension of Good's idea.  
Also, conjectures for the corresponding $q$-analogs
are supplied.  Finally, perturbed versions of the $q$-Dixon summation
formula are presented.
\end{abstract}
\end{frontmatter}

\section{Introduction}

\subsection{Notation}
For $n$ a nonnegative integer, we define the following symbols:
  \begin{align}
     &{\av} := \langle a_1, a_2, \dots, a_n \rangle, 
         \tag{$n$-vector of symbolic nonnegative integers} \\
     &{\xv} := \langle x_1, x_2, \dots, x_n \rangle, 
          \tag{$n$-vector of indeterminants}\\
     &{\zerov} := \langle 0,0, \dots, 0 \rangle,
          \tag{$n$-dimensional zero vector}\\
     &{\ev}_{k}:= \langle 0,0,\dots,0,1,0,0,\dots,0 \rangle,
          \tag{the $n$-vector with 1 in the $k$th position and 0 elsewhere}\\
     &\s := a_1 + a_2 + \cdots + a_n, 
     \tag{first elementary symmetric polynomial in $n$ indeterminants}\\
     &(A;q)_n := \prod_{i=0}^{n-1} (1-Aq^i), \tag{rising $q$-factorial}\\
     &F_n(\xv ; \av) := \prod_{1\leqq i<j \leqq n} 
   \left(1-\frac{x_i}{x_j} \right)^{a_j}  \left(1-\frac{x_j}{x_i}\right)^{a_i}
       \tag{Dyson product}\\
     &\mathcal{F}_n( \xv ; \av ; q):= 
        \prod_{1\leqq i<j \leqq n} 
        \left( \frac{x_i q}{x_j};q \right)_{a_j}  
        \left(\frac{x_j}{x_i};q \right)_{a_i}, \label{qDysonProd}
  \tag{$q$-Dyson product}
   \end{align}
and let  $[Y]Z $ denote the coefficient of $Y$ in the expression $Z$, thus e.g.
\begin{gather*}
   [x^3 y^2] (3 + 5 x^3 y^2 - 6 xy) = 5,\\
  [1] (3 + 5 x^3 y^2 - 6 xy) = [x^0 y^0] (3 + 5 x^3 y^2 - 6 xy) = 3,\\
  [xy^2](3+5 x^3 y^2 - 6xy) = 0.
\end{gather*}

\subsection{Background}
F.~J.~Dyson~\cite[p. 152, Conjecture C]{fjd} 
conjectured that the constant term in the 
Laurent polynomial 
$\prod_{1\leqq i<j \leqq n} 
       \left(1-\frac{x_i}{x_j} \right)^{a_j}  \left(1-\frac{x_j}{x_i}\right)^{a_i}$
is the multinomial coefficient; i.e.

\begin{DC}
  For $n \in \Zset_+$, 
    \begin{equation} \label {dc}
        [1] F_n(\xv; \av) = 
        \frac{\s!}{a_1! a_2! \cdots a_n!}.
     \end{equation}
\end{DC}
Dyson's conjecture~\eqref{dc} was first proved independently by 
J.~Gunson~\cite{jg}
and K.~Wilson~\cite{kw}.  
Later I.~J.~Good~\cite{ijg} supplied the most compact and elegant proof.

G.~E.~Andrews~\cite[p. 216]{gea} extended~\eqref{dc} to a $q$-analog:

\begin{qDC}
 For $n \in \Zset_+$, 
  \begin{equation} \label {qdc}
        [1] \mathcal{F}_n(\xv; \av ; q) = 
        \frac{(q;q)_{\s}}{(q;q)_{a_1} (q;q)_{a_2} \cdots (q;q)_{a_n} }.
     \end{equation}
\end{qDC}

  The first proof of~\eqref{qdc} was given by
D.~Zeilberger and D.~M.~Bressoud~\cite{zb}.  Recently, another proof
was given by I.~M.~Gessel and G.~Xin~\cite{gx}.

  In~\cite{sz}, together with Zeilberger, I showed that with the aid of our 
Maple/Mathematica packages
\texttt{GoodDyson}, the computer can, subject only to 
limitations of time and memory capacity, conjecture a closed form expression
for

 \[ [x_1^{b_1} x_2^{b_2} \cdots x_n^{b_n}] F_n(\xv; \av), \]
 and \textit{automatically supply a proof} for any \emph{fixed} 
positive integer $n$ and
\emph{fixed} vector $\bv = \langle b_1, b_2, \dots, b_n \rangle$.

\subsection{Theorems and Conjectures}
 The results of~\cite{sz} are extended here to \emph{generic} $n$ for certain
 vectors $\bv$, and a corresponding $q$-analog is conjectured for each.
I made heavy use of Maple in forming these conjectures.  
I will prove
 \begin{thm}\label{1m1} Let $r$ and $s$ be fixed integers with 
$1\leqq r\neq s \leqq n$ and $n\geqq 2$.  Then
 \begin{equation}\label{thm1}
 [x_r/x_s] F_n(\xv ; \av) = -\left(\frac{a_s}{1+\s -a_s}\right)
 \frac{\s!}{a_1! a_2!\cdots a_n!}.
\end{equation} 
 \end{thm}
and provide a conjecture for its $q$-analog:
 \begin{conj}[$q$-analog of Theorem~\ref{1m1}]\label{q1m1}
 Let $r$ and $s$ be fixed integers with $1\leqq r\neq s \leqq n$
and $n\geqq 2$.  Then
 \[ [x_r/x_s] \mathcal{F}_n(\xv ; \av ; q) = 
 -q^{L(r,s)}\left( \frac{1-q^{a_s}}{1-q^{1+\s-a_s}}\right)
 \frac{(q;q)_{\s}}{(q;q)_{a_1} (q;q)_{a_2}\cdots (q;q)_{a_n} }, \]
 where 
   \[ L(r,s) = \left\{ \begin{array}{ll}
                      1+\s-\sum_{k=r}^s a_k, &\mbox{if $r<s$}\\
                       \sum_{k=s+1}^{r-1} a_k, &\mbox{if $r>s$.}
                    \end{array} \right.\]
 \end{conj}
  
\begin{rem}
Notice that the right hand side of Eq.~\eqref{thm1} is \emph{independent
of $r$,} the subscript of the variable which appears to a positive power.  
In other words, $[x_k/x_s] F_n(\xv; \av)$ is the same for all
$k\neq s$. 
This can be explained by the fact that the only factors contributing to
the $x_k/x_s$ term in the expansion of $F_n(\xv; \av)$ are
 \[\underset{i\neq k}{\prod_{i=1}^n } \left( 1- \frac{x_i}{x_s} \right)^{a_s}, 
\]
which is clearly invariant under any permutation of the subscripts of
the $x_i$.  The analogous phenomenon occurs in Theorems~\ref{2m1m1} and
\ref{11m1m1} as well.
\end{rem}

Next, we have 
 \begin{thm}\label{2m1m1}
  Let $r$, $s$, and $t$ be distinct fixed integers with 
$1\leqq r, s, t \leqq n$ and $n\geqq 3$.  Then
 \[ \left[ \frac{x_r^2}{x_s x_t} \right] F_n(\xv ; \av) = 
 \left(\frac{a_s a_t \Big( (1+\s) + (1+\s-a_s-a_t) \Big) }
 {(1+\s-a_s-a_t)(1+\s-a_s)(1+\s-a_t)}\right)
 \frac{\s!}{a_1! a_2!\cdots a_n!}, \]
 \end{thm}
and the following conjecture for its $q$-analog: 
\begin{conj}[$q$-analog of Theorem~\ref{2m1m1}]\label{q2m1m1}
  Let $r$, $s$, and $t$ be distinct fixed integers with
$1\leqq r, s, t \leqq n$ and $n\geqq 3$.  Without loss of
generality we may assume that $s<t$.  Then
 \begin{gather*}
\left[ \frac{x_r^2}{x_s x_t} \right] \mathcal{F}_n(\xv ; \av ; q) =
 q^{L(r,s,t)} 
 \left(\frac{(1-q^{a_s})(1-q^{a_t}) \Big( (1-q^{1+\s}) + 
q^{M(r,s,t)}(1-q^{1+\s-a_s-a_t}) 
\Big) }
 {(1-q^{1+\s-a_s-a_t})(1-q^{1+\s-a_s})(1-q^{1+\s-a_t}) }\right) \\ \times
 \frac{(q;q)_{\s}}{(q;q)_{a_1} (q;q)_{a_2}\cdots (q;q)_{a_n}},
\end{gather*} 
where
   \[ L(r,s,t) = \left\{ \begin{array}{ll}
          2+2\s-2\sum_{k=r}^t a_k+\sum_{k=s+1}^{t-1}a_k, &\mbox{if $r<s<t$,}\\
          1+\s-\sum_{k=s}^t a_k + 2\sum_{k=s+1}^{r-1}a_k,&\mbox{if $s<r<t$,}\\   
          2\sum_{k=t+1}^{r-1} a_k+\sum_{k=s+1}^{t-1} a_k,&\mbox{if $s<t<r$,}
                    \end{array} \right.\] and
   \[ M(r,s,t) = \left\{ \begin{array}{ll}
            a_t, &\mbox{if $r<s<t$ or $s<t<r$,}\\
            a_s, &\mbox{if $s<r<t$.}
                      \end{array}\right.\]
 \end{conj}

Finally, we have
\begin{thm}\label{11m1m1}
  Let $r$, $s$, $t$, and $u$ be distinct fixed integers with 
$1\leqq r, s, t, u \leqq n$ and $n\geqq 4$.
  Then
 \[ \left[ \frac{x_r x_s}{x_t x_u} \right] F_n(\xv ; \av) =
 \left(\frac{a_t a_u \Big( (1+\s) + (1+\s-a_t-a_u) \Big) }
 {(1+\s-a_t-a_u)(1+\s-a_t)(1+\s-a_u)}\right)
 \frac{\s!}{a_1! a_2!\cdots a_n!}. \]
 \end{thm}

\begin{conj}[$q$-analog of Theorem~\ref{11m1m1}]\label{q11m1m1}
  Let $r$, $s$, $t$ and $u$ be distinct fixed integers with
$1\leqq r, s, t, u \leqq n$ and $n\geqq 4$.  Without loss of
generality we may assume that $r<s$ and $t<u$.  Then
 \begin{gather*}
\left[ \frac{x_r x_s}{x_t x_u} \right] \mathcal{F}_n(\xv ; \av ; q) =
 q^{L(r,s,t,u)}
 \left(\frac{(1-q^{a_t})(1-q^{a_u}) \Big( (1-q^{1+\s}) +
q^{M(r,s,t,u)}(1-q^{1+\s-a_t-a_u})
\Big) }
 {(1-q^{1+\s-a_t-a_u})(1-q^{1+\s-a_t})(1-q^{1+\s-a_u}) }\right) \\ \times
 \frac{(q;q)_{\s}}{(q;q)_{a_1} (q;q)_{a_2}\cdots (q;q)_{a_n}},
\end{gather*}
where
   \[ L(r,s,t,u) = \left\{ \begin{array}{ll}
     2+2\s-2\sum_{k=r}^u a_k+\sum_{k=r}^{s-1}a_k+\sum_{k=t+1}^{u-1} a_k, 
                        &\mbox{if $r<s<t<u$,}\\
          1+\s-\sum_{k=r}^{u} a_k+\sum_{k=t+1}^{s-1}a_k , 
                      &\mbox{if $r<t<s<u$,}\\
          1+\s-\sum_{k=r}^{s-1} a_k + 2\sum_{k=t+1}^{r-1} a_k + 
               \sum_{k=t+1}^{u-1} a_k +2\sum_{k=u+1}^{s-1} a_k ,
                       &\mbox{if $r<t<u<s$,}\\
          1+\s-\sum_{k=t}^{u} a_k+ \sum_{k=r}^{s-1} a_k 
           +2\sum_{k=t+1}^{r-1} a_k,
                    &\mbox{if $t<r<s<u$,}\\
          \sum_{k=t+1}^{r-1} a_k + \sum_{k=u+1}^{s-1}a_k, &\mbox{if $t<r<u<s$,}\\
          \sum_{k=r}^{s-1}a_k+\sum_{k=t+1}^{u-1} a_k+ 2\sum_{k=u+1}^{r-1}a_k,
                        &\mbox{if $t<u<r<s$,}
                    \end{array} \right.\] and
   \[ M(r,s,t,u) = \left\{ \begin{array}{ll}
            a_u,         &\mbox{if $r<s<t<u$ or $r<t<u<s$ or $t<u<r<s$,}\\
            1+\s  &\mbox{if $r<t<s<u$ or $t<r<u<s$,}\\
            a_t,  &\mbox{if $t<r<s<u$.}\\
                      \end{array}\right.\]
 \end{conj}

\begin{rem} Certain special cases of Conjectures~\ref{q1m1},~\ref{q2m1m1},
and~\ref{q11m1m1} have been proved by John 
Stembridge~\cite[p. 347, Cor. 7.4]{jrs}.
Stembridge proved that in the case where
$\av = \langle a, a, \dots, a \rangle$, and
$b_{\rho+1} = b_{\rho+2} = \cdots = b_{\rho+\tau} = -1$, for
$\rho$ and $\tau$ satisfying 
$0\leqq \rho\leqq n$ and $1\leqq\tau\leqq n-\rho$, 
\begin{equation} \label{stemb}
[x_1^{b_1} x_2^{b_2} \cdots x_n^{b_n}]\mathcal{F}_n(\xv; \av; q)
= (-1)^\tau q^{b_1+b_2+\cdots+b_\rho + am} 
\frac{(q;q)_{an} (q^a; q^a)_\tau (q;q^a)_{\rho+\sigma} }
{ (q;q)_a^n (q;q^a)_n },
\end{equation}     
where $m = \sigma\tau + \sum_{i=1}^\rho (i-1) b_i - \sum_{i=1}^{n-\rho-\tau}
i\ b_{n-i+1}$.
Conjectures~\ref{q1m1},~\ref{q2m1m1},
and~\ref{q11m1m1} do indeed agree with~\eqref{stemb} where they overlap,
which, of course, provides some evidence in favor of the conjectures.
\end{rem}

  The theorems will be proved in \S\ref{GoodProofs}.
Special cases of the conjectured $q$-analogs will be discussed in 
some detail in \S\ref{qDixon}, followed by some concluding remarks
in \S\ref{conc}.
 
 \section{Generalized Good Proofs}\label{GoodProofs}
 \subsection{Good's proof of Dyson's conjecture}
 It will be instructive to review the proof of~\eqref{dc} due 
to Good~\cite{ijg}
presented in a way that will make it easy
to see how it naturally generalizes to the variations of Dyson's conjecture
under consideration here.  The proof divides neatly into three parts: recurrence,
initial condition, and boundary conditions.
Let \[ c_n^\bv(\av):=  [x_1^{b_1} x_2^{b_2} \cdots x_n^{b_n}] F_n(\xv; \av).\]
Thus Dyson's conjecture is the assertion that
 \[ c_n^\zerov(\av)= \frac{\s!}{a_1! a_2! \cdots a_n!}.\]
    
 \subsubsection{Recurrence}
For $a_1, a_2, \dots, a_n > 0$, we have, by Lagrange
interpolation,
     \begin{equation}\label{Frec}
    F_n(\xv; \av) = \sum_{k=1}^n
    F_n(\xv; \av - \ev_k).
      \end{equation}
    Thus the same recurrence must hold term by term when~\eqref{Frec} is 
 expanded, and in particular the recurrence must hold for the constant term, so
 we have
  \begin{equation}\label{crec}\tag{$R$}
     c_n^{\zerov}(\av) = \sum_{k=1}^n
     c_n^{\zerov}(\av - \ev_k).
  \end{equation}
  
    \subsubsection{Initial Condition}  It is easily verified that
      \begin{equation}\label{ic}\tag{$I$}
     c_n^{\zerov}(\zerov) = 1.
       \end{equation}
      
     \subsubsection{Boundary Conditions}
     For $k$ fixed and $1\leqq k \leqq n$, 
      \begin{gather}
    F_n(\xv; \langle a_1, a_2, \dots, a_{k-1},0,a_{k+1}, \dots a_n \rangle ) 
    \nonumber\\
     = F_{n-1}(\langle x_1, x_2, \dots, x_{k-1},x_{k+1},\dots,x_n\rangle ;
    \langle a_1,\dots,a_{k-1},a_{k+1}, \dots a_n\rangle )
     \left\{ \underset{i\neq k}{ \prod_{i=1}^n} 
     \frac{(x_i-x_k)^{a_i}}{x_i^{a_i}}\right\}
     \label{bc1}
     \end{gather}
     Notice that we have segregated 
the factors involving $x_k$ (those in braces)
from those which are independent of $x_k$. Find the
Taylor expansion of $  \underset{i\neq k}{ \prod_{i=1}^n} 
\frac{(x_i-x_k)^{a_i}}{x_i^{a_i}}$ about $x_k=0$. Extract the
  coefficient of $x_k^0$ from both sides of~\eqref{bc1} to obtain
      \begin{gather}
    [x_k^0]F_n(\xv; \langle a_1, a_2, \dots, a_{k-1},0,a_{k+1}, \dots a_n \rangle ) 
    \nonumber\\
     = P_k^\zerov \times 
   F_{n-1}(\langle x_1, x_2, \dots, x_{k-1},x_{k+1},\dots,x_n\rangle ;
    \langle a_1,\dots,a_{k-1},a_{k+1}, \dots a_n\rangle )
     \label{bc11},
     \end{gather}  
 where \begin{equation}\label{Pkdef}
     P_k^\bv = [x_k^{b_k}]  \underset{i\neq k}{ \prod_{i=1}^n} 
  \frac{(x_i-x_k)^{a_i}}{x_i^{a_i}}.
   \end{equation}
In the case of Dyson's original conjecture, we have
  $P_k^\zerov = 1$ for all $k$ and $n$.
 
   Apply the constant term operator to both sides of~\eqref{bc11} to obtain
   \begin{equation}\label{bc}\tag{$B$}
   c_n^\zerov(\langle a_1, a_2, \dots, a_{k-1},0, a_{k+1}, \dots, a_n \rangle)
  =  c_{n-1}^\zerov(\langle a_1, a_2, \dots, a_{k-1},a_{k+1},\dots a_n \rangle)  
   \end{equation}
 for $k=1,2,\dots, n$.

   Finally, since~\eqref{crec},~\eqref{ic}, and~\eqref{bc} uniquely
 determine $c_n^\zerov(\av),$ and the multinomial coefficient
 $\s!/ a_1!\cdots a_n!$ also satisfies ~\eqref{crec},~\eqref{ic}, and~\eqref{bc},
the result follows. \qed
 
 \subsection{Proof of Theorem~\ref{1m1}}
 Theorem~\ref{1m1} asserts that if $\bv=\ev_r - \ev_s$,
 \begin{equation}\label{c1m1}
 c_n^{\bv}(\av) = -\left(\frac{a_s}{1+\s-a_s}\right)
 \frac{\s!}{a_1! a_2!\cdots a_n!}.
\end{equation}
 
  \subsubsection{Recurrence}
        It was already noted that by Lagrange interpolation, 
for $a_1, a_2, \dots, a_n > 0$, we have
     \begin{equation}\label{Frec2}
    F_n(\xv; \av) = \sum_{k=1}^n
    F_n(\xv; \av - \ev_k).
      \end{equation}
    Thus the same recurrence must hold term by term when~\eqref{Frec2} is 
 expanded, and in particular the recurrence must hold for the $x_r/x_s$ term,
and so
  \begin{equation}\label{crec2}\tag{$R'$}
     c_n^{\ev_r - \ev_s}(\av) = \sum_{k=1}^n
    c_n^{\ev_r - \ev_s}(\av-\ev_k).
  \end{equation}
  
    \subsubsection{Initial Condition}
      \begin{equation}\label{ic2}\tag{$I'$}
     c_n^{\ev_r - \ev_s}(\zerov) = 0.
       \end{equation}
   
       \subsubsection{Boundary Conditions}
     For $k$ fixed and $1\leqq k \leqq n$,
      \begin{gather}
    F_n(\xv; \langle a_1, a_2, \dots, a_{k-1},0,a_{k+1}, \dots a_n \rangle ) 
    \nonumber\\
     = F_{n-1}(\langle x_1, x_2, \dots, x_{k-1},x_{k+1},\dots,x_n\rangle ;
    \langle a_1,\dots,a_{k-1},a_{k+1}, \dots a_n\rangle )
     \left\{ \underset{i\neq k}{ \prod_{i=1}^n} 
     \frac{(x_i-x_k)^{a_i}}{x_i^{a_i}}\right\}
     \label{bc21}
     \end{gather}
  Once again, we have segregated the factors involving $x_k$ (those in braces)
from those which are independent of $x_k$. Next, find the
Taylor expansion of $  \underset{i\neq k}{ \prod_{i=1}^n} 
\frac{(x_i-x_k)^{a_i}}{x_i^{a_i}}$ about $x_k=0$. Extract the
coefficient of $x_k^{b_k}$ from both sides of~\eqref{bc21} to obtain
      \begin{gather}
    [x_k^{b_k}]F_n(\xv; \langle a_1, a_2, \dots, a_{k-1},0,a_{k+1}, \dots a_n \rangle ) 
    \nonumber\\
     = P_k^{\bv} \times 
   F_{n-1}(\langle x_1, x_2, \dots, x_{k-1},x_{k+1},\dots,x_n\rangle ;
    \langle a_1,\dots,a_{k-1},a_{k+1}, \dots a_n\rangle )
     \label{bc20}
     \end{gather}  
 where
  \[ P_k^{\bv} = \left\{
                          \begin{array}{ll}
                           -\underset{i\neq k}{\sum_{i=1}^n} \frac{a_i}{x_i},
                               &\mbox{if $k=r$,}\\
                            0, &\mbox{if $k=s$,}\\
                            1, &\mbox{otherwise,} 
                       \end{array} \right.  \]
    and thus by extracting the coefficient of $x_r x_s^{-1} x_k^{b_k}$ 
from both sides of~\eqref{bc20},
we obtain

\begin{multline}
    c_n^{\ev_r - \ev_s} 
      (\langle a_1, a_2, \dots, a_{k-1}, 0 , a_{k+1}, \dots a_n \rangle)\\
 =   \left\{
         \begin{array}{ll}
         -\underset{i\neq k}{\sum_{i=1}^n} a_i c_{n-1}^{\ev_i^{(k)}-\ev_s^{(k)}}
         (\langle a_1, a_2, \dots, a_{k-1} , a_{k+1}, \dots a_n \rangle), 
            &\mbox{if $k=r$,}\\
         0, &\mbox{if $k=s$,}\\
         c_{n-1}^{\ev_r^{(k)} - \ev_s^{(k)}} 
             (\langle a_1, a_2, \dots, a_{k-1} , a_{k+1}, \dots a_n \rangle),
         &\mbox{otherwise,}
  \end{array} \right.
\label{bc2}\tag{$B'$}
\end{multline}
where \[ \ev_j^{(k)} = \langle \delta_{1,j}, \delta_{2,j}, \dots,
  \delta_{k-1,j}, \delta_{k+1,j}, \dots, \delta_{n,j} \rangle, \] 
with $\delta_{i,j}$ denoting the Kronecker delta function.

\subsubsection{ The RHS of~\eqref{c1m1} also satisfies (R), (I), and (B)}
  Since ~\eqref{crec2},~\eqref{ic2}, and~\eqref{bc2} uniquely
 determine $c_n^{\ev_r - \ev_s}(\av),$ once we establish that 
 $d^{\ev_r-\ev_s}_n(\av):=-\left(\frac{a_s}{1+\s-a_s}\right)
 \left( \frac{\s!}{a_1! a_2!\cdots a_n!}\right)$ also satisfies ~\eqref{crec2},~\eqref{ic2}, 
 and~\eqref{bc2}, the result will follow.   
While this fact may not be obious \emph{a priori}, we shall soon see
that nothing beyond elementary algebra is required to establish its
truth.

 Without loss of generality, we may assume that $r=1$ and $s=n$, for if
not, the indeterminants in $F_n(\xv;\av)$ may be relabeled accordingly.
We note that
\[ d_n^{\ev_1 - \ev_n} (\av) = -\left( \frac{a_n}{1+a_1+a_2+\dots+a_{n-1}} 
\right) \left(\frac{ \sigma_n(\av)!}{a_1! a_2! \cdots a_{n}!} \right) . \]   
\begin{align*}
 \sum_{k=1}^n d^{\ev_1-\ev_n}_n (\av-\ev_k) &= 
-\frac{(a_n-1) (a_1 + \cdots +a_n -1)!}{(1+a_1+\cdots a_{n-1}) a_1!
\cdots a_{n-1}! (a_n - 1)!}\\
& \qquad - \sum_{k=1}^{n-1} \frac{a_k a_n (a_1+\cdots+a_n-1)!}
{(a_1+\cdots+a_{n-1}) a_1!\cdots a_n!}\\
&=\frac{-a_n (a_1+\cdots+a_n-1)!}{(1+a_1+\cdots+a_{n-1})a_1!\cdots a_n!
(a_1+\cdots+a_{n-1})}\\
& \qquad \times\Big\{ (a_n-1)(a_1+\cdots+a_{n-1})+ 
  \sum_{k=1}^{n-1} a_k (1+a_1+\cdots+a_{n-1}) \Big\}\\
&=\frac{-a_n (a_1+\cdots+a_n-1)!}{(1+a_1+\cdots+a_{n-1})a_1!\cdots a_n!
(a_1+\cdots+a_{n-1})}\\
& \qquad \times \Big\{ (a_1+\cdots+a_{n-1})(a_n-1 + 1+a_1 +\cdots a_{n-1})
\Big\}\\
&= \frac{ -a_n (a_1+\cdots a_n)!}{ (1+a_1+\cdots + a_{n-1}) a_1! \cdots a_n!}
\\
&= d_n^{\ev_1 - \ev_n}(\av),
\end{align*}  and thus ($R'$) is satisfied.

 Clearly,
\[ d_n^{\ev_1 - \ev_n}(\zerov) = 0, \] so ($I'$) is satisfied.

 Also,
\begin{align*}
&\quad-\sum_{i=2}^n a_i d_{n-1}^{\ev_i^{(1)}-\ev_n^{(1)}} (\langle a_2, \dots, a_n
\rangle)\\
 &= -a_n d^{\zerov}_{n-1} (\langle a_2, \dots, a_n \rangle)
-\sum_{i=2}^{n-1} a_i d^{\ev_i^{(1)}-\ev_n^{(1)}} (\langle a_2, \dots, a_n
\rangle)\\
&= \frac{(a_2 + \cdots a_n)!}{a_2! \cdots a_n!}
\left( \frac{a_2 a_n}{1+a_2+\cdots a_{n-1}} + \cdots +\frac{a_{n-1}a_n}
{1+a_2+\cdots + a_{n-1}} - a_n \right)\\
&=\frac{(a_2 + \cdots a_n)! a_n}{a_2! \cdots a_n! (1+a_2+\cdots+a_{n-1}} 
\big( a_2 + \cdots + a_{n-1} - (1 + a_2 + \cdots + a_{n-1}) \big)\\
&= -\frac{ (a_2 + \cdots + a_n )! a_n}{ a_2! \cdots a_n! (1+a_2 + \cdots
+ a_{n-1} )}\\
&=d_n^{\ev_1^{(1)} - \ev_n^{(1)}} (\langle 0, a_2, \dots, a_n \rangle),
\end{align*}
and thus $d_n^{\ev_r - \ev_s}(\av)$ satisfies $(B')$ when $a_r = 0$.

 Clearly,
\[ d_n^{\ev_1^{(n)}-\ev_n^{(n)}} (\langle a_1, \dots, a_{n-1}, 0 \rangle)
= 0, \] and so $d_n^{\ev_r - \ev_s}(\av)$ satisfies $(B')$ when $a_s = 0$.
 
 Finally, for $1<k<n$, we have
\begin{align*}
& =\quad d_n^{\ev_1^{(k)} - \ev_n^{(k)}} (\langle a_1, \dots, a_{k-1}, 0,
a_{k+1} , \dots, a_n \rangle) \\
&= \frac{ -a_n}{1+a_1+\cdots + a_{k-1} + a_{k+1} + \cdots a_n}
\frac{ (a_1 + \cdots + a_{k-1} + a_{k+1} + \cdots + a_n)!}{
a_1! \cdots a_{k-1}! a_{k+1}! \cdots a_n!}\\
&= d_{n-1}^{\ev_1^{(k)} - \ev_n^{(k)}}(\langle a_1, \dots, a_{k-1},
a_{k+1}, \dots, a_n \rangle ) ,
\end{align*}
where $d_n^\zerov (\av) = \s!/a_1!\cdots a_n!$ by~\eqref{dc},
and thus $d_n^{\ev_r - \ev_s}(\av)$ satisfies $(B')$ when $k$ is 
different from both $r$ and $s$. \qed
 
\begin{rem}\label{remark}
Clearly, the only nontrivial difference between the proof of~\eqref{dc}
and that of Theorem~\ref{1m1} lies in the observation that $P_k^{\bv}$
(see~\eqref{Pkdef}) varies with $\bv$.  Once $P_k^{\bv}$ is known for
a given $\bv$, the boundary condition (\eqref{bc} and \eqref{bc2}
in the two previous cases) follows immediately.   
\end{rem}

\subsection{Proof of Theorem~\ref{2m1m1}}
In light of Remark~\ref{remark}, we need only supply $P_k^{\bv}$, for
$\bv = 2\ev_r - \ev_s - \ev_t$. 

\[ P_k^{2\ev_r-\ev_s-\ev_t} = 
            \left\{
               \begin{array}{ll}
                  \left(\underset{i\neq k}{\sum_{i=1}^n} 
                     \frac{a_i(a_i-1)}{2 x_i^2}
               + \underset{i\neq k}{\sum_{1\leqq i<j\leqq n}}
                      \frac{a_i a_j}{x_i x_j} \right),
                               &\mbox{if $k=r$,}\\
                            0, &\mbox{if $k=s$ or $k=t$,}\\
                            1, &\mbox{otherwise,}
                       \end{array} \right.  \]
which implies

\begin{multline*}
    c_n^{2\ev_r-\ev_s-\ev_t} 
        (\langle a_1, a_2, \dots, a_{k-1}, 0 , a_{k+1}, \dots a_n \rangle)
 \\=   \left\{
         \begin{array}{ll}
         \underset{i\neq k}{\sum_{i=1}^n} \frac{a_i(a_i-1)}{2}
              c_{n-1}^{2\ev_i^{(k)}i-\ev_s^{(k)} - \ev_t^{(k)} }
         (\langle a_1, a_2, \dots, a_{k-1} , a_{k+1}, \dots a_n \rangle)\\
       \qquad +\underset{i\neq k}{\sum_{1\leqq i<j\leqq n}} 
              a_i a_j c_{n-1}^{\ev_i^{(k)}+\ev_j^{(k)}-\ev_s^{(k)}-\ev_t^{(k)}}
               (\langle a_1, a_2, \dots, a_{k-1} , a_{k+1}, \dots a_n \rangle),
            &\mbox{if $k=r$,}\\
         0, &\mbox{if $k=s$ or $k=t$,}\\
         c_{n-1}^{2\ev_r^{(k)}-\ev_s^{(k)}-\ev_t^{(k)}}
             (\langle a_1, a_2, \dots, a_{k-1} , a_{k+1}, \dots a_n \rangle),
         &\mbox{otherwise.}
  \end{array} \right. 
  \end{multline*}

\subsection{Proof of Theorem~\ref{11m1m1}}
Similarly,
\[ P_k^{\ev_r+\ev_s-\ev_t-\ev_u} =
            \left\{
               \begin{array}{ll}
                  \left(-\underset{i\neq k}{\sum_{i=1}^n}
                     \frac{a_i}{x_i}
                      \right),
                               &\mbox{if $k=r$ or $k=s$,}\\
                            0, &\mbox{if $k=t$ or $k=u$,}\\
                            1, &\mbox{otherwise,}
                       \end{array} \right.  \]
which implies

\begin{multline*}
    c_n^{\ev_r+\ev_s-\ev_t-\ev_u}
        (\langle a_1, a_2, \dots, a_{k-1}, 0 , a_{k+1}, \dots a_n \rangle)
\\
 =   \left\{
         \begin{array}{ll}
         \underset{i\neq k}{-\sum_{i=1}^n} a_i
              c_{n-1}^{\ev_s^{(k)}+\ev_i^{(k)}-\ev_t^{(k)}-\ev_u^{(k)}}
         (\langle a_1, a_2, \dots, a_{k-1} , a_{k+1}, \dots a_n \rangle)
            &\mbox{if $k=r$, }\\
         \underset{i\neq k}{-\sum_{i=1}^n} a_i
              c_{n-1}^{\ev_r^{(k)}+\ev_i^{(k)}-\ev_t^{(k)}-\ev_u^{(k)}}
         (\langle a_1, a_2, \dots, a_{k-1} , a_{k+1}, \dots a_n \rangle)
            &\mbox{if $k=s$, }\\
         0, &\mbox{if $k=t$ or $k=u$,}\\
         c_{n-1}^{\ev_r^{(k)}+\ev_s^{(k)}-\ev_t^{(k)}-\ev_u^{(k)}}
             (\langle a_1, a_2, \dots, a_{k-1} , a_{k+1}, \dots a_n \rangle),
         &\mbox{otherwise.}
  \end{array} \right. 
  \end{multline*}

\section{Perturbed versions of $q$-Dixon}\label{qDixon}
It is well known (see~\cite{gea}) that the $n=3$ case of the $q$-Dyson
conjecture is equivalent to a $q$-analog of a hypergeometric
summation formula of A.~C.~Dixon~\cite{acd}. 

This is because 
 \begin{align*}
& \mathcal{F}_3 (\langle x,y,z \rangle; \langle a,b,c \rangle) \\
&=  (y/x;q)_a (z/x;q)_a (xq/y;q)_b (z/y;q)_b (xq/z;q)_c (yq/z;q)_c \\
&= \frac{(-1)^{b+2c} q^{\binom{b}{2}+2\binom{c}{2}}} {x^{2a} y^{2b} z^{2c}}
    \prod_{i=0}^{a+b-1} (x-yq^{i-b})
    \prod_{i=0}^{a+c-1} (x-zq^{i-c})
    \prod_{i=0}^{b+c-1} (y-zq^{i-c})\\
&= \sum_{h,i,j \geqq 0} \gp{a+b}{h}{q} \gp{a+c}{i}{q} \gp{b+c}{j}{q}
  (-1)^{b+2c+h+i+j} q^{\binom{b-h}{2}+ \binom{c-i}{2} + \binom{c-j}{2}}
  x^{b+c-h-i} y^{-b+c+h-i} z^{-2c+i+j}, 
 \end{align*}
where the last equality follows from a triple application of a
corollary of the $q$-binomial theorem due to Rothe (see
\cite[p. 490, Cor. 10.2.2 (c)]{aar}), and 
\[ \gp{A}{B}{q} =  \left\{
   \begin{array}{ll}
     \frac{(q;q)_A}{(q;q)_B (q;q)_{A-B}} &\mbox{if $0\leqq A\leqq B$}\\
     0 &\mbox{otherwise.}
   \end{array} \right.
\]
It is then a straightforward exercise in linear algebra combined with
the change of variable $k=j+c$ to obtain

\begin{multline*}
 \left[ \frac{x^\alpha y^\beta}{z^{\alpha+\beta}} \right] 
    \mathcal{F}_3 (\langle x,y,z \rangle; \langle a,b,c \rangle; q)
 \\ = \sum_{k\in\Zset} \gp{a+b}{k+b+\beta}{q}
                         \gp{b+c}{k+c}{q}
                         \gp{c+a}{k+a+\alpha+\beta}{q}
           (-1)^{k+\alpha} 
     q^{\binom{k+1}{2} +  \binom{k+1+\beta}{2} + \binom{k+\alpha+\beta}{2} }. 
\end{multline*}
For $\alpha=\beta=0$, combined with the $n=3$ case of the $q$-Dyson
theorem, we obtain the $q$-Dixon sum of
Andrews~\cite[p. 216, equation (5.6)]{gea},
which he proved using the $q$-Pfaff-Saalsch\"utz summation
(see \cite[equation (II.12)]{gr}.)

  Similarly, the following six identities follow from the $n=3$ case
of Conjecture~\ref{q1m1}:
\begin{gather} \label{first}
\sum_{k\in\Zset} \gp{a+b}{k+b-1}{q}
                         \gp{b+c}{k+c}{q}
                         \gp{c+a}{k+a}{q}
           (-1)^{k} q^{k(3k-1)/2}
= \gp{a+b+c}{a,b,c}{q} \left( \frac{1-q^b}{1-q^{1+a+c}}\right) q^{1+c}\\
\sum_{k\in\Zset} \gp{a+b}{k+b}{q}
                         \gp{b+c}{k+c}{q}
                         \gp{c+a}{k+a+1}{q}
           (-1)^{k} q^{3k(k+1)/2-1}
= \gp{a+b+c}{a,b,c}{q} \left( \frac{1-q^c}{1-q^{1+a+b}}\right) \\
\sum_{k\in\Zset} \gp{a+b}{k+b+1}{q}
                         \gp{b+c}{k+c}{q}
                         \gp{c+a}{k+a}{q}
           (-1)^{k} q^{3k(k+1)/2 + 1}
=\gp{a+b+c}{a,b,c}{q} \left( \frac{1-q^a}{1-q^{1+b+c}}\right)\\
\sum_{k\in\Zset} \gp{a+b}{k+b}{q}
                         \gp{b+c}{k+c}{q}
                         \gp{c+a}{k+a-1}{q}
           (-1)^{k} q^{k(3k-1)/2+1}
= \gp{a+b+c}{a,b,c}{q} \left( \frac{1-q^a}{1-q^{1+b+c}}\right) q^{b}\\
\sum_{k\in\Zset} \gp{a+b}{k+b+1}{q}
                         \gp{b+c}{k+c}{q}
                         \gp{c+a}{k+a+1}{q}
           (-1)^{k+1} q^{k(3k+5)/2 }
= \gp{a+b+c}{a,b,c}{q} \left( \frac{1-q^c}{1-q^{1+a+b}}\right) q^{a}\\
\sum_{k\in\Zset} \gp{a+b}{k+b-1}{q}
                         \gp{b+c}{k+c}{q}
                         \gp{c+a}{k+a-1}{q}
           (-1)^{k+1} q^{3k(k-1)/2 + 1}
= \gp{a+b+c}{a,b,c}{q} \left( \frac{1-q^b}{1-q^{1+a+c}}\right),      
\end{gather}
where $$\gp{a+b+c}{a,b,c}{q} = \frac{(q;q)_{a+b+c}}{(q;q)_a (q;q)_b (q;q)_c}.$$

  The corresponding 
identities arising from the $n=3$ case of Conjecture~\ref{q2m1m1}
are
\begin{multline}
\sum_{k\in\Zset} \gp{a+b}{k+b-1}{q}
                         \gp{b+c}{k+c}{q}
                         \gp{c+a}{k+a+1}{q}
           (-1)^{k} q^{3k(k-1)/2}
\\= \gp{a+b+c}{a,b,c}{q} 
   \frac{(1-q^b)(1-q^c)}{(1-q^{1+b})(1-q^{1+a+b})(1-q^{1+a+c})}
    \left( (1-q^{1+a+b+c}) - q^c(1-q^a) \right) 
\end{multline} 

\begin{multline}
\sum_{k\in\Zset} \gp{a+b}{k+b+2}{q}
                         \gp{b+c}{k+c}{q}
                         \gp{c+a}{k+a+1}{q}
           (-1)^{k+1} q^{k(3k+7)/2 + 2}
\\= \gp{a+b+c}{a,b,c}{q}
   \frac{(1-q^a)(1-q^c)}{(1-q^{1+b})(1-q^{1+a+b})(1-q^{1+b+c})}
    \left( (1-q^{1+a+b+c}) - q^a(1-q^b) \right) 
\end{multline}

\begin{multline}
\sum_{k\in\Zset} \gp{a+b}{k+b-1}{q}
                         \gp{b+c}{k+c}{q}
                         \gp{c+a}{k+a-2}{q}
           (-1)^{k+1} q^{k(3k-5)/2 + 3}
\\= \gp{a+b+c}{a,b,c}{q}
   \frac{(1-q^a)(1-q^b)}{(1-q^{1+c})(1-q^{1+a+c})(1-q^{1+b+c})}
    \left( (1-q^{1+a+b+c}) - q^b(1-q^c) \right) \label{last}
\end{multline}

\begin{rem}
    Each of the identities \eqref{first} through \eqref{last}
is a ${}_3\phi_2$ summation formula, and as such is automatically
verifiable by the $q$-WZ algorithm of Wilf and Zeilberger~\cite{wz:multi}.  
It is well known
that Zeilberger's algorithm and its $q$-analog does not always find the minimal
order recurrence satisfied by a given summand (see, e.g. ~\cite{gea:pfaff3}
or~\cite[p. 116 ff.]{pwz}).   In each case considered here,
the $q$-Zeilberger algorithm, as implemented in Maple by Zeilberger's
package \texttt{qEKHAD} and in \emph{Mathematica} by A.~Riese's package
\texttt{qZeil.m} (see \cite{pr}), 
a recurrence of order at least three was found for
the sum side, even though there \emph{must} be a first order recurrence
since the right hand side is a sum of a fixed number of finite products.  
Even Paule's creative
symmetrization technique (see~\cite[section 5.2]{pr}) does not improve the  
order of the recurrence in these examples.
\end{rem}

\begin{rem}
The same technique could be used to produce $q$-hypergeometric summation
formulas corresponding to the case $n=4$.  Here the resulting sum sides
would be triple sums, and one could attempt to obain automated proofs
of these in \emph{Mathematica} using Riese's\texttt{qMultiSum.m} 
package of~\cite{ar},
or in Maple using Zeilberger's~\texttt{qMultiZeilberger} package
~\cite{dz:qMZ}.        

  Due to computer memory and time limitations, it is highly doubtful that 
the identities corresponding to $n>4$ could be successfully handled
on today's computers.
\end{rem}

\section{Conclusion}\label{conc}
The obvious next step is to try to find proofs for the conjectured
$q$-analogs.  A combinatorial proof would be particularly nice,
since would potentially explain the role played by the factors
$q^L$ and $q^M$ in the conjectures, a feature that disappears in
the ordinary $q=1$ case.

\begin{ack}
I thank Doron Zeilberger for getting me interested in the $q$-Dyson
conjecture, and for the many stimulating conversations which occured
as a result.  I also with to express my gratitude to the referees, whose
thoughtful suggestions helped to improve the paper.
\end{ack}

\end{document}